\documentclass[11pt]{article}

\usepackage[T2A]{fontenc}
\usepackage[cp1251]{inputenc}
\usepackage{amsfonts}
\usepackage{eufrak}
\usepackage{amssymb}
\usepackage{amsmath}

\begin{document}

\newtheorem{theorem}{Theorem}
\newtheorem{lemma}{Lemma}
\newtheorem{corollary}{Corollary}
\newtheorem{definition}{Definition}
\newtheorem{proposition}{Proposition}
\newtheorem{remark}{Remark}

\begin{center}
	\textbf{\Large IDENTICALLY DISTRIBUTED RANDOM VECTORS ON LOCALLY COMPACT ABELIAN GROUPS}
	
	\bigskip
	
	{\Large Margaryta Myronyuk}
	
	\bigskip
	
	\textit{B. Verkin Institute for Low Temperature Physics and Engineering\\ of the National Academy of
		Sciences of Ukraine, \\ Nauky Ave. 47, Kharkiv, Ukraine}
	
	\bigskip
	
	\textit{Bielefeld University, \\ Universit\"{a}tsstra\ss{}e 25, Bielefeld, Germany}
	
	\bigskip
	
	myronyuk@ilt.kharkov.ua
\end{center}

\begin{abstract}
	L. Klebanov proved the following theorem. Let $\xi_1, \dots, \xi_n$ be independent random variables. Consider linear forms $L_1=a_1\xi_1+\cdots+a_n\xi_n,$ 
	$L_2=b_1\xi_1+\cdots+b_n\xi_n,$
	$L_3=c_1\xi_1+\cdots+c_n\xi_n,$
	$L_4=d_1\xi_1+\cdots+d_n\xi_n,$  
	where the coefficients $a_j, b_j, c_j, d_j$ are real numbers. If the random vectors $(L_1,L_2)$ and $(L_3,L_4)$ are identically distributed, then all $\xi_i$ for which $a_id_j-b_ic_j\neq 0$ for all $j=\overline{1,n}$ are Gaussian random variables. The present article is devoted to an
	analog of the Klebanov theorem in the case when random variables take
	values in a locally compact Abelian group and the coefficients of
	the linear forms are integers.
\end{abstract}



\bigskip

\section{Introduction}

Many researches are devoted to the studying of linear forms of independent real-valued random variables (\cite{KaLiRa}). The analytic theory of them started in the famous paper by Cramer. He proved the following decomposition theorem for a Gaussian distribution (\cite{Cramer}).

\medskip

\noindent \textbf{The Cramer theorem.} \textit{If a linear form $L=\xi_1+\xi_2$ has a Gaussian distribution then the independent random variables $\xi_1, \xi_2$ also have Gaussian distributions. }

\medskip

\noindent Kac (\cite{Kac}) and Bernstein (\cite{Bernstein}) considered the sum $L_1=\xi_1+\xi_2$ and the difference $L_2=\xi_1-\xi_2$ and proved that the independence of the linear forms $L_1$ and $L_2$ implies that the independent random variables $\xi_1, \xi_2$ have Gaussian distributions. Darmois (\cite{Darmois}) and Skitovich (\cite{Skitovich}) generalized the Kac-Bernstein theorem for two arbitrary linear forms.

\medskip

\noindent \textbf{The Darmois-Skitovich theorem.} \textit{ Let $a_jb_j\neq 0$ for all $j=\overline{1,n}$. If linear forms $L_1=a_1\xi_1+\cdots+a_n\xi_n$ and $L_2=b_1\xi_1+\cdots+b_n\xi_n$ are independent then the independent random variables $\xi_1,\dots, \xi_n$ have Gaussian distributions. }

\medskip

\noindent Heyde  considered another statistical property of two arbitrary linear forms and proved the following theorem (\cite{Heyde}).

\medskip

\noindent \textbf{The Heyde theorem.} \textit{ Let $a_jb_j\neq 0$ and $a_ib_j+a_jb_i\neq 0$ for all $i,j=\overline{1,n}$. If the conditional distribution of $L_2=b_1\xi_1+\cdots+b_n\xi_n$ given $L_1=a_1\xi_1+\cdots+a_n\xi_n$ is symmetric then the independent random variables $\xi_1,\dots, \xi_n$ have Gaussian distributions. }

\medskip

\noindent Klebanov  considered four arbitrary linear forms of independent random variables and obtained the following result  (\cite{Klebanov}).

\medskip

\noindent \textbf{The Klebanov theorem.} \textit{ Let $\xi_1, \dots, \xi_n$ be independent random variables. Consider linear forms $L_1=a_1\xi_1+\cdots+a_n\xi_n,$ 
	$L_2=b_1\xi_1+\cdots+b_n\xi_n,$
	$L_3=c_1\xi_1+\cdots+c_n\xi_n,$
	$L_4=d_1\xi_1+\cdots+d_n\xi_n,$  
	where the coefficients $a_j, b_j, c_j, d_j$ are real numbers. If the random vectors $(L_1,L_2)$ and $(L_3,L_4)$ are identically distributed, then all $\xi_i$ for which $a_id_j-b_ic_j\neq 0$ for all $j=\overline{1,n}$ are Gaussian random variables. }

\medskip

\noindent Note that the Klebanov theorem implies the Darmois-Skitovich theorem and the Heyde theorem. These and other researches devoted to the studying of linear forms of independent real-valued random variables can be founded e.g. in \cite{Bernstein}-\cite{Darmois}, \cite{Heyde}, \cite{Ibragimov}, \cite{KS}-\cite{LukacsKing}, \cite{Skitovich}, see also \cite{KaLiRa} for more references.

Many branches of the classical probability theory and mathematical statistics have been generalized on different algebraic structures.
In particular, there exist a lot of papers devoted to the studying of linear forms of independent random variables with values in locally compact Abelian groups. The group analogue or the Cramer theorem was proved in the paper \cite{FeKr}. Many papers are devoted to the generalizations of the Kac-Bernstein theorem, the Darmois-Skitovich theorem and the Heyde theorem (see e.g. \cite{FeKr}, \cite{Fe1990}, \cite{FeGr2010}, \cite{My2006}, \cite{Myronyuk2020}, see also \cite{FeBook2} and \cite{FeBook3} for more references).   

The aim of the present paper is to 
 prove the analog of the Klebanov theorem in the case
when random variables take values in a locally compact Abelian group
and coefficients of the linear forms are integers.

\section{Notation and definitions}

In the article we use standard results on abstract harmonic analysis
(see e.g. \cite{HeRo1}).

Let $X$ be a second countable locally compact  Abelian group,
$Y=X^\ast$ be its character group, and  $(x,y)$ be the value of a
character $y \in Y$ at an element $x \in X$. Denote by $Aut(X)$ the group of all topological automorphisms of the group $X$. Let $H$ be a subgroup
of $Y$. Denote by $A(X,H)=\{x \in X: (x,y)=1 \ \ \forall \ y \in
H\}$ the annihilator of $H$. For each integer $n$, $n \ne 0,$ let
$f_n : X \mapsto X$ be the endomorphism $f_n x=nx.$ Set $X^{(n)} =
f_n(X)$, $X_{(n)}=Ker f_n$. Denote by ${\mathbb{Z}}(n)$  the finite
cyclic group of order $n$. 

Let ${M^1}(X)$ be the convolution semigroup of probability
distributions on $X$,  $$\widehat \mu(y) = \int_X (x, y) d\mu(x)$$ be
the characteristic function of a distribution $\mu \in {M^1}(X)$,
and $\sigma(\mu)$ be the support of $\mu$. If $H$ is a closed
subgroup of $Y$ and $\widehat \mu(y)=1$ for $y \in H$, then
$\widehat\mu(y+h) = \widehat\mu(y)$ for all $y \in Y$, $ h \in H$ and
$\sigma(\mu) \subset A(X, H)$. For $\mu \in {M^1}(X)$ we define the
distribution $\bar \mu \in M^1(X)$ by the rule $\bar \mu(B) =
\mu(-B)$ for all Borel sets $B \subset X$. Note that $\widehat
{\bar \mu}(y) = \overline{\widehat \mu(y)}$.

Let $x\in X$. Denote by $E_x$ the degenerate distribution
concentrated at the point $x$, and by $D(X)$ the set of all
degenerate distributions on
$X$. A distribution $\gamma \in {M^1}(X)$ is called Gaussian (\cite[\S
4.6]{Pa}) if its characteristic function can be represented in the
form
\begin{equation}\label{i1}
	\widehat\gamma(y)= (x,y)\exp\{-\varphi(y)\},
\end{equation}	
where $x \in X$ and $\varphi(y)$ is a continuous nonnegative
function satisfying the equation
\begin{equation}\label{i2}
\varphi(u+v)+\varphi(u-v)=2[\varphi(u)+
\varphi(v)], \quad u, \ v \in Y.
\end{equation}
Denote by $\Gamma(X)$ the set
of Gaussian distributions on $X$. We note that according to this
definition $D(X)\subset \Gamma(X)$. Denote by $I(X)$ the set of
shifts of Haar distributions $m_K$ of compact subgroups $K$ of the
group $X$. Note that

\begin{equation}
\widehat{m}_K(y)=
\left\{%
\begin{array}{ll}
1, & \hbox{$y\in A(Y,K)$;} \\
0, & \hbox{$y\not\in A(Y,K)$.} \\
\end{array}%
\right.
\end{equation}
We note that if a distribution $\mu \in \Gamma(X)*I(X)$, i.e.
$\mu=\gamma*m_K$, where $\gamma\in \Gamma(X)$, then $\mu$ is
invariant with respect to the compact subgroup $K \subset X$ and
under the natural homomorphism $X \mapsto X/K \ \mu$ induces a
Gaussian distribution on the factor group $X/K$. Therefore the class
$\Gamma(X)*I(X)$ can be considered as a natural analogue of the
class $\Gamma(X)$ on locally compact Abelian groups. 

Let $f(y)$ be a function on $Y$, and $h\in Y.$ Denote by $\Delta_h$
the finite difference operator
$$\Delta_h f(y)=f(y+h)-f(y).$$
A function $f(y)$ on $Y$ is called a polynomial if
$$\Delta_{h}^{l+1}f(y)=0$$
for some $n$ and for all $y,h \in Y$. The minimal $l$ for which  this equality holds is called the degree of the polynomial $f(y)$.

An integer $a$ is said to be admissible for a group $X$ if $X^{(a)}
\ne \{0\}$. The admissibility of integers $a$
is a group analogue of the condition  $a\ne 0$ for the case of  $X = \mathbb{R}$.

\section{Main results}

Let $X$ be a second countable locally compact Abelian group. Let $\xi_1, \dots, \xi_n$ be independent random variables with values in $X$. Consider the linear forms 
\begin{equation}\label{th1.1.1}
L_1=a_1\xi_1+\cdots+a_n\xi_n,
\end{equation}
\begin{equation}\label{th1.1.2}
L_2=b_1\xi_1+\cdots+b_n\xi_n,
\end{equation}
\begin{equation}\label{th1.1.3}
L_3=c_1\xi_1+\cdots+c_n\xi_n,
\end{equation}
\begin{equation}\label{th1.1.4}
L_4=d_1\xi_1+\cdots+d_n\xi_n,
\end{equation}
where the coefficients $a_j, b_j, c_j, d_j$ are 
integers. We will keep the designation $L_1$, $L_2$, $L_3$ and $L_4$ throughout the article. 

The following theorem gives the description of locally compact Abelian groups for which the group analog of the Klebanov theorem takes place.

\begin{theorem}\label{theorem1}
	Let $X$ be a second countable locally
	compact Abelian group. Let $\xi_1, \dots, \xi_n$ be independent random variables with values in $X$ and distributions $\mu_{\xi_j}$ with non-vanishing characteristic functions. Consider the linear forms $(\ref{th1.1.1})-(\ref{th1.1.4})$, 
	where the coefficients $a_j, b_j, c_j, d_j$ are 
	integers. If the random vectors $(L_1,L_2)$ and $(L_3,L_4)$ are identically distributed, then for all those $i$ for which
	\begin{equation}\label{th1.1}
	a_i d_j - b_i c_j \text{ are admissible integers for } X \text{ for all } j=\overline{1,n},
	\end{equation}
	the following
	statements hold:
	
	\qquad $(i)$ If $X$ is a torsion-free group then 
	$\mu_{\xi_i}\in\Gamma(X)$;
	
	\qquad $(ii)$ If $X=X_{(p)}$ then $\mu_{\xi_i}\in D(X)$.	
\end{theorem}


\begin{remark}
	{\rm The class of torsion-free locally compact Abelian groups is very wide. For example, it includes the real line $\mathbb{R}$, the group of integers   $\mathbb{Z}$, the group of rational numbers $\mathbb{Q}$, compact torsion-free Abelian groups, their products. Note that each compact torsion free Abelian group is topologically isomorphic to a group of the form }
	\begin{equation*}
	\Sigma_{\boldmath{a}}^\mathfrak{m} \times \textbf{P}_{p\in \textit{P}} \Delta_p^{\mathfrak{m}_p},
	\end{equation*}
	{\rm where $\textit{P}$ is the set of prime numbers, $\mathfrak{m}$ and $\mathfrak{m}_p$ are arbitrary cardinal
	numbers, $\Sigma_{\boldmath{a}}$ is the $\boldmath{a}$-adic solenoid ($\boldmath{a}=(2,3,4,\dots)$) and $\Delta_p$ is the group of $p$-adic integers $(\cite[\S 25.8]{HeRo1})$.}
\end{remark}

\begin{remark}
	{\rm If $X= X_{(p)}$, where $p$ is a fixed prime number, then $X$ is topologically isomorphic to
	the group
	\begin{equation}\label{rem1}
	\mathbb{Z}(p)^\mathfrak{m} \times \mathbb{Z}(p)^{\mathfrak{n}*},
	\end{equation}
	where $\mathfrak{m}$ and $\mathfrak{n}$ are arbitrary cardinal
	numbers, ${\mathbb{Z}(p)}^{\mathfrak{m}}$ is considered in the
	product topology, and $ {\mathbb{Z}(p)}^{\mathfrak{n}*}$ is
	considered in the discrete topology $(\cite[\S 25.29]{HeRo1})$.}
\end{remark}

We need some lemmas. The following lemma allows us to reduce the proof of Theorem \ref{theorem1} to the study of solutions of some functional equation.

\begin{lemma}\label{Equation}
	Let $X$ be a second countable locally
	compact Abelian group, $Y=X^*$. Let $\xi_1, \dots, \xi_n$ be independent random variables with values in $X$ and distributions $\mu_{\xi_j}$. Consider the linear forms $(\ref{th1.1.1})-(\ref{th1.1.4})$. The random vectors $(L_1,L_2)$ and $(L_3,L_4)$ are identically distributed if and only if the characteristic functions $\widehat\mu_{\xi_j}(y)$ satisfy the equation
	\begin{equation}\label{th1.2}
	\prod_{j=1}^{n}  \widehat\mu_{\xi_j} (a_j u+ b_j v)=\prod_{j=1}^{n}  \widehat\mu_{\xi_j} (c_j u+ d_j v),\quad u,v\in Y.
	\end{equation}
\end{lemma}

\textit{\textbf{Proof.}} Using the independence of $\xi_1, ..., \xi_n$, we obtain from the definition that
\begin{eqnarray*}
	\widehat\mu_{(L_1,L_2)}(u,v)=\mathbf{E}\left[ (L_1,L_2) (u,v) \right] =\mathbf{E}\left[ (L_1,u) (L_2,v) \right] = \\ =\mathbf{E}\left[ (a_1\xi_1+\cdots+a_n\xi_n,u) (b_1\xi_1+\cdots+b_n\xi_n,v) \right] = \\ =\mathbf{E}\left[ (a_1\xi_1,u)\cdots(a_n\xi_n,u) (b_1\xi_1,v)\cdots(b_n\xi_n,v) \right] = \\ =\mathbf{E}\left[ (\xi_1,a_1 u)\cdots(\xi_n, a_n u) (\xi_1,b_1 v)\cdots(\xi_n, b_n v) \right]= \\ =\mathbf{E}\left[ (\xi_1,a_1 u +b_1 v)\cdots(\xi_n, a_n u+b_n v) \right]= \prod_{j=1}^{n}  \widehat\mu_{\xi_j} (a_j u+ b_j v),\quad u,v\in Y.
\end{eqnarray*}
Analogously we obtain that
\begin{equation*}
\widehat\mu_{(L_3,L_4)}(u,v)=\prod_{j=1}^{n}  \widehat\mu_{\xi_j} (c_j u+ d_j v),\quad u,v\in Y.
\end{equation*} 
Since the vectors $(L_1,L_2)$ and $(L_3,L_4)$ are identically distributed, we have $\widehat\mu_{(L_1,L_2)}(u,v)=\widehat\mu_{(L_3,L_4)}(u,v)$. Thus, equation
(\ref{th1.2}) is valid.
\hfill $\square$

\medskip

The following lemma is a group analogue of the Cramer theorem on the decomposition of a Gaussian  distribution.

\begin{lemma}\label{GrKramer} $(\cite{FeKr})$ Let $X$ be a second countable locally compact Abelian group containing no subgroup topologically isomorphic to the circle group. Let $\gamma\in \Gamma(X)$ and $\gamma=\gamma_1*\gamma_2$, where $\gamma_1, \gamma_2 \in M^1(X)$. Then $\gamma_1, \gamma_2 \in \Gamma(X)$.
\end{lemma}

The following lemma is a group analogue of the Marcinkiewicz theorem.

\begin{lemma}\label{GrMarc} $(\cite{FeMarc})$
	Let $X$ be a second countable locally compact Abelian group containing no subgroup topologically isomorphic to the circle group. Let $\mu \in M^1(X)$ and the characteristic function $\widehat{\mu}(y)$ is of the form
	\begin{equation*}
		\widehat{\mu}(y)=\exp\{ \varphi(y)\}, \quad y\in Y,
	\end{equation*}
	where $\varphi(y)$ is a continuous polynomial. Then $\mu\in \Gamma(X)$.
\end{lemma}

\begin{lemma}\label{FinDifEq} 
	Let $Y$ be an Abelian group. Let $\varphi_{j}(y)$ and $\psi_{j}(y)$ be functions on $Y$ which satisfy the equation 
	\begin{equation}\label{fde1}
	\sum_{j=1}^{m}  \varphi_{j} (a_j u+ b_j v)=\sum_{j=1}^{n}  \psi_{j} (c_j u+ d_j v)+q(u,v),\quad u,v\in Y,
	\end{equation}
	where the coefficients $a_j, b_j, c_j, d_j$ are integers and $q(u,v)$ is a continuous polynomial of degree $l$ on the group
	$Y^2$. 
	Then each function $\varphi_j(y)$ satisfies the equation
	\begin{equation}\label{fde2}
	\Delta_{a_jh+b_jk}^{l+1}\Delta_{(a_jb_{i_1}-b_ja_{i_1})l_{i_1}} ... \Delta_{(a_jb_{i_{m-1}}-b_ja_{i_{m-1}})l_{i_{m-1}}} \Delta_{(a_jd_{1}-b_jc_{1})k_{1}} ...  \Delta_{(a_jd_n-b_jc_n)k_n} \varphi_{j} (a_j u+b_jv)=$$$$=0,\quad u,v\in Y,
	\end{equation}
	where $h, k, k_1, \dots, k_n$ and $l_{i_1}, \dots, l_{i_{m-1}}$ are arbitrary elements of $Y$ and indexes $i_1, \dots, i_{m-1}$ take values $1,\dots,j-1,j+1,\dots,m$. Moreover, if $q(u,v)\equiv 0$ then operator $\Delta_{a_jh+b_jk}^{l+1}$ is missing from formula $(\ref{fde2})$.
	 
\end{lemma}

\textbf{\textit{Proof.}} Let $k_n$ be an arbitrary element of $Y$. Substitute $u+d_n k_n$ for $u$ and $v-c_n k_n$ for $v$ in (\ref{fde1}) and subtract (\ref{fde1}) from the resulting equation. We obtain
\begin{equation}\label{fde3}
\sum_{j=1}^{m} \Delta_{(a_jd_n-b_jc_n)k_n} \varphi_{j} (a_j u+ b_j v)=\sum_{j=1}^{n-1} \Delta_{(c_jd_n-d_jc_n)k_n} \psi_{j} (c_j u+ d_j v)+ \Delta_{(d_n k_n, -c_n k_n)} q(u,v),\quad u,v\in Y.
\end{equation}
The right-hand side of equation (\ref{fde3}) no longer contains the function $\psi_n$.

Let $k_{n-1}$ be an arbitrary element of $Y$. Substitute $u+d_{n-1}k_{n-1}$ for $u$ and $v-c_{n-1}k_{n-1}$ for $v$ in (\ref{fde3}) and subtract (\ref{fde3}) from the resulting equation. We obtain
\begin{equation}\label{fde4}
\sum_{j=1}^{m} \Delta_{(a_jd_{n-1}-b_jc_{n-1})k_{n-1}} \Delta_{(a_jd_n-b_jc_n)k_n} \varphi_{j} (a_j u+ b_j v)=$$$$ \sum_{j=1}^{n-2} \Delta_{(c_jd_{n-1}-d_jc_{n-1})k_n} \Delta_{(c_jd_n-d_jc_n)k_n} \psi_{j} (c_j u+ d_j v) +$$$$+ \Delta_{(d_{n-1}k_{n-1}, -c_{n-1}k_{n-1})}\Delta_{(d_n k_n, -c_n k_n)} q(u,v),\quad u,v\in Y.
\end{equation}
The right-hand side of equation (\ref{fde4}) no longer contains the function $\psi_{n-1}$.

After $n$ steps we come to the equation of the form
\begin{equation}\label{fde5}
\sum_{j=1}^{m} \Delta_{(a_jd_{1}-b_jc_{1})k_{1}} ... \Delta_{(a_jd_{n-1}-b_jc_{n-1})k_{n-1}} \Delta_{(a_jd_n-b_jc_n)k_n} \varphi_{j} (a_j u+ b_j v)=$$$$=\Delta_{(d_{1}k_{1}, -c_{1}k_{1})}...\Delta_{(d_n k_n, -c_n k_n)} q(u,v),\quad u,v\in Y,
\end{equation}
where $k_1, \dots, k_n$ are arbitrary elements of $Y$.

Let $l_m$ be an arbitrary element of $Y$. Substitute $u+b_ml_m$ for $u$ and $v-a_ml_m$ for $v$ in (\ref{fde5}) and subtract (\ref{fde5}) from the resulting equation. We obtain
\begin{equation}\label{fde6}
\sum_{j=1}^{m-1} \Delta_{(a_jb_m-b_ja_m)l_m} \Delta_{(a_jd_{1}-b_jc_{1})k_{1}} ...  \Delta_{(a_jd_n-b_jc_n)k_n} \varphi_{j} (a_j u+ b_j v)=$$$$= \Delta_{(b_ml_m, -a_ml_m)} \Delta_{(d_{1}k_{1}, -c_{1}k_{1})}...\Delta_{(d_n k_n, -c_n k_n)} q(u,v),\quad u,v\in Y.
\end{equation}
The left-hand side of equation (\ref{fde6}) no longer contains the function $\varphi_m$.

After $n+m-1$ steps we come to the equation of the form
\begin{equation}\label{fde7}
\Delta_{(a_1b_{2}-b_1a_{2})l_{2}} ... \Delta_{(a_1b_m-b_1a_m)l_m} \Delta_{(a_1d_{1}-b_1c_{1})k_{1}} ...  \Delta_{(a_1d_n-b_1c_n)k_n} \varphi_{1} (a_1 u+ b_1 v)=$$$$= \Delta_{(b_2l_2, -a_2l_2)} ... \Delta_{(b_ml_m, -a_ml_m)} \Delta_{(d_{1}k_{1}, -c_{1}k_{1})}...\Delta_{(d_n k_n, -c_n k_n)} q(u,v),\quad u,v\in Y,
\end{equation}
where $k_1, \dots, k_n$ and $l_2, \dots, l_m$ are arbitrary elements of $Y$. 

Since $q(u,u)$ is a polynomial of degree $l$, applying the operator
$\Delta^{l+1}_{(h, k)}$ to (\ref{fde7}), we get equation (\ref{fde2}) for the function $\varphi_{1}(y)$.   By analogy we prove that each function $\varphi_j(y)$ satisfies equation (\ref{fde2}).
\hfill $\square$

\medskip

\textbf{\textit{Proof of Theorem \ref{theorem1}.}} 
It follows from Lemma \ref{Equation} that the characteristic functions $\widehat\mu_{\xi_j}(y)$ satisfy equation (\ref{th1.2}). Renumbering functions in (\ref{th1.2}), we can assume that condition (\ref{th1.1}) is fulfilled for all $i=\overline{1,m}$ and it is not fulfilled for all $i=\overline{m+1,n}$, $m\leq n$.

\textit{Case 1.} First we consider case 1 when there exist admissible integers $b_ia_j - b_j a_i$ for $i\neq j$ ($i,j=\overline{1,m}$).



(i) Let $X$ be a torsion-free group. In this case if an integer $b_i
a_j - b_j a_i$ is not admissible for $X$ then $b_i a_j - b_j a_i=0$. Also if an integer $a$ is admissible for $X$ then 

\begin{equation}\label{th1.10.-1}
\overline{Y^{(a)}}=Y \text{ for all nonzero integeres } a.
\end{equation}

\noindent Note that if there exists $j_0$ ($j_0\in \{1,\dots,n\}$) such that $c_{j_0}=d_{j_0}=0$ then there does not exist such $i_0$ ($i_0\in \{1,\dots,m\}$) that condition (\ref{th1.1}) is fulfilled. Therefore we can suppose from the beginning that there do not exist $j_0$ ($j_0\in \{1,\dots,n\}$) such that $c_{j_0}=d_{j_0}=0$.


Renumbering functions in (\ref{th1.2}), we can assume that
\begin{equation*}
	a_1,\dots,a_t \neq 0, \quad a_{t+1}=\dots=a_m=0,
\end{equation*}
\begin{equation*}
b_1,\dots,b_s \neq 0, \quad b_{s+1}=\dots=b_t=0, \quad b_{t+1},\dots, b_m \neq 0,
\end{equation*}
where $s\leq t \leq m$.

If $a_{t+1}=\dots=a_m=0$, then condition (\ref{th1.1}) implies that

\begin{equation} \label{th1.10.1}
	c_j\neq 0 \text{ for all } j=\overline{1,n}. 
\end{equation}

If $b_{s+1}=\dots=b_t=0$, then condition (\ref{th1.1}) implies that

\begin{equation} \label{th1.10.2}
d_j\neq 0 \text{ for all } j=\overline{1,n}. 
\end{equation}

Also renumbering first $s$ functions in (\ref{th1.2}), we can assume that

\begin{equation}\label{th1.11}
{a_1 \over b_1}= \cdots = {a_{r_1} \over b_{r_1}}=\alpha_1, {a_{r_1+1} \over b_{r_1+1}}= \cdots = {a_{r_2} \over b_{r_2}}=\alpha_2,\cdots,
{a_{r_{k}+1} \over b_{r_{k}+1}}= \cdots = {a_s \over b_s}=\alpha_{k+1},
\end{equation}
where $1\leq r_1<r_2<\cdots < r_k< r_{k+1} \leq s$, and $\alpha_i\neq \alpha_j$ for all $i
\ne j$.

Consider integers

\begin{equation}\label{th1.12}
A={b_1 b_{r_1+1}\cdots b_{r_k+1}}, \quad B={a_1 a_{r_1+1}\cdots a_{r_k+1}},
\end{equation}

\begin{equation}\label{th1.13}
A_0={A\over b_1}, \quad A_j={A\over b_{r_j+1}}, \quad j=1,2,\dots, k,
\end{equation} 

\begin{equation}\label{th1.14}
B_0={B\over a_1}, \quad B_j={B\over a_{r_j+1}}, \quad j=1,2,\dots, k,
\end{equation} 

Substitute $A u$ for $u$ and $B v$ for $v$ in (\ref{th1.2}). We get

\begin{equation}\label{th1.14.1}
\prod_{j=1}^{s}  \widehat\mu_{\xi_j} (a_jA u+ b_jB v) \prod_{j=s+1}^{t}  \widehat\mu_{\xi_j} (a_jA u) \prod_{j=t+1}^{m}  \widehat\mu_{\xi_j} (b_jB v)\prod_{j=m+1}^{n}  \widehat\mu_{\xi_j} (a_jA u+ b_jB v) =$$$$=\prod_{j=1}^{n}  \widehat\mu_{\xi_j} (c_j Au+ d_j Bv),\quad u,v\in Y.
\end{equation}
Taking into account (\ref{th1.12})-(\ref{th1.14}), we rewrite (\ref{th1.14.1}) in the following form

\begin{equation}\label{th1.15}
\prod_{j=1}^{r_1}  \widehat\mu_{\xi_j} (a_jb_1A_0 u+ b_ja_1B_0 v) \prod_{j=r_1+1}^{r_2}  \widehat\mu_{\xi_j} (a_jb_{r_1+1}A_1 u+ b_ja_{r_1+1}B_1 v) \cdots $$$$ \cdots \prod_{j=r_k+1}^{s}  \widehat\mu_{\xi_j} (a_jb_{r_k+1}A_k u+ b_ja_{r_k+1}B_k v) \prod_{j=s+1}^{t}  \widehat\mu_{\xi_j} (a_jA u) \prod_{j=t+1}^{m}  \widehat\mu_{\xi_j} (b_jB v) \times  $$$$ \times\prod_{j=m+1}^{n}  \widehat\mu_{\xi_j} (a_jA u+ b_jB v) =\prod_{j=1}^{n}  \widehat\mu_{\xi_j} (C_j u+ D_j v),\quad u,v\in Y,
\end{equation}
where $C_j=c_j A$, $D_j= d_j B$. Taking into account (\ref{th1.11}), we obtain from (\ref{th1.15})

\begin{equation}\label{th1.16}
\prod_{j=1}^{r_1}  \widehat\mu_{\xi_j} \left(a_1b_j(A_0 u+B_0 v)\right) \prod_{j=r_1+1}^{r_2}  \widehat\mu_{\xi_j} \left(a_{r_1+1} b_j(A_1 u+ B_1 v)\right) \cdots $$$$ \cdots \prod_{j=r_k+1}^{s}  \widehat\mu_{\xi_j} \left(a_{r_k+1}b_{j}(A_k u+ B_k v)\right) \prod_{j=s+1}^{t}  \widehat\mu_{\xi_j} (a_jA u) \prod_{j=t+1}^{m}  \widehat\mu_{\xi_j} (b_jB v)\times  $$$$\times \prod_{j=m+1}^{n}  \widehat\mu_{\xi_j} (a_jA u+ b_jB v) =\prod_{j=1}^{n}  \widehat\mu_{\xi_j} (C_j u+ D_j v),\quad u,v\in Y.
\end{equation}
Put 
$$\widehat\mu_0(y)= \prod_{j=1}^{r_1}  \widehat\mu_{\xi_j} \left(a_1b_jy\right), \quad \widehat\mu_1(y)= \prod_{j=r_1+1}^{r_2}  \widehat\mu_{\xi_j} \left(a_{r_1+1} b_jy \right), \cdots ,\quad \widehat\mu_k(y)= \prod_{j=r_k+1}^{s}  \widehat\mu_{\xi_j} \left(a_{r_k+1}b_{j}y\right),$$

$$ \widehat\mu_{k+1}(y) = \prod_{j=s+1}^{t}  \widehat\mu_{\xi_j} (a_jA y), \quad \widehat\mu_{k+2}(y) = \prod_{j=t+1}^{m}  \widehat\mu_{\xi_j} (b_jB y).  $$ 

\noindent Note that each $\widehat\mu_i(y)$ is a characteristic function of a distribution $\mu_{i}$ on $X$. Namely,

	$$\mu_{0}=f_{a_1b_1}(\mu_{\xi_1})*\cdots*f_{a_1b_{r_1}}(\mu_{\xi_{r_1}}),$$$$ \mu_{1}=f_{a_{r_1+1}b_{r_1+1}}(\mu_{\xi_1})*\cdots* f_{a_{r_1+1}b_{r_2}}(\mu_{\xi_{r_1}}),$$$$ \cdots,$$$$ \mu_{k}=f_{a_{r_k+1}b_{r_k+1}}(\mu_{\xi_{r_k+1}})*\cdots* f_{a_{r_k+1}b_{s}}(\mu_{\xi_{s}}), $$$$
	\mu_{k+1}=f_{a_{s+1}A}(\mu_{\xi_{s+1}})*\cdots* f_{a_{t}A}(\mu_{\xi_{t}}),
	$$$$
	\mu_{k+2}=f_{b_{t+1}B}(\mu_{\xi_{t+1}})*\cdots* f_{b_{m}B}(\mu_{\xi_{m}}).$$
If we prove that all $\widehat\mu_{i}(y)$ ($i=0,\overline{1,k+2}$) are characteristic functions of Gaussian distributions then, taking into account Lemma \ref{GrKramer} and ({\ref{th1.10.-1}), we obtain that all distributions  $\mu_{\xi_j} \in \Gamma(X)$ ($j=\overline{1,m}$).  Now we get from (\ref{th1.16})

\begin{equation}\label{th1.17}
\prod_{i=0}^{k}  \widehat\mu_{i} \left(A_i u+B_i v\right) \ \widehat\mu_{k+1} (A_{k+1}u+B_{k+1}v)   \widehat\mu_{k+2} (A_{k+2}u+B_{k+2}v)  \prod_{j=m+1}^{n}  \widehat\mu_{\xi_j} (A_j u+ B_j v)= $$$$ =\prod_{j=1}^{n}  \widehat\mu_{\xi_j} (C_j u+ D_j v),\quad u,v\in Y,
\end{equation}
where $A_{k+1}=1$, $B_{k+1}=0$; $A_{k+2}=0$, $B_{k+2}=1$; $A_j=a_j A$, $B_j=b_j B$ ($j=\overline{m+1,n}$).
We have that 

\begin{equation}\label{th1.18}
A_iD_j-B_iC_j \neq 0, \quad i=\overline{0,k+2} \text{ and } j=\overline{1,n}.
\end{equation}
Indeed, taking into account (\ref{th1.1}), (\ref{th1.10.1}), (\ref{th1.10.2}), (\ref{th1.12})-(\ref{th1.14}), we get
\begin{equation*}
A_iD_j-B_iC_j = { A\over b_{r_i+1}} \cdot d_j B - { B\over a_{r_i+1} }\cdot c_j A= {AB\over a_{r_i+1}b_{r_i+1}} (a_{r_i+1} d_j - b_{r_i+1}  c_j )  \neq 0,\quad i=\overline{0,k}, $$$$
A_{k+1}D_j-B_{k+1}C_j = D_j=d_j B \neq 0,$$$$
A_{k+2}D_j-B_{k+2}C_j = -C_j=- c_j A  \neq 0.
\end{equation*}
Also we have that 
\begin{equation}\label{th1.19}
A_iB_j-B_iA_j \neq 0, \quad i\neq j, \quad i=\overline{0,k+2}, \quad j=\overline{0,k+2} \text{ and } j=\overline{m+1,n}.
\end{equation}
Indeed, taking into account (\ref{th1.11}), we have
\begin{equation*}
A_iB_j-B_iA_j = { A\over b_{r_i+1}} \cdot {B\over a_{r_j+1}} - {B\over a_{r_i+1}} \cdot {A\over b_{r_j+1} }= {AB ( a_{r_i+1}  b_{r_j+1} - b_{r_i+1} a_{r_j+1}) \over b_{r_i+1} a_{r_j+1} a_{r_i+1} b_{r_j+1}}  \neq 0, \ i,j=\overline{0,k}, i\neq j, $$$$
A_iB_{k+1}-B_iA_{k+1} = -B_i \neq 0, \ i=\overline{0,k}, $$$$
A_iB_{k+2}-B_iA_{k+2} = A_i  \neq 0, \quad i=\overline{0,k}, $$$$
A_{k+1}B_{k+2}-B_{k+1}A_{k+2} = A_{k+1}B_{k+2}=1  \neq 0.
\end{equation*}
We make some notes about $a_j, b_j$ when $j=\overline{m+1,n}$. First we note that if $a_{j_0}=b_{j_0}=0$ for some $j_0$ then the corresponding multiplier $\widehat\mu_{\xi_{j_0}}$ is absent in equation (\ref{th1.17}). Therefore we can ignore this case. Since condition (\ref{th1.1}) is not valid for such $a_j, b_j$, we have $a_jd_{\alpha}-b_jc_{\alpha}=0$ for some $c_{\alpha}, d_{\alpha}$. Since $c_{\alpha}, d_{\alpha}$ are not zero simultaneously, we can suppose for definiteness that $d_{\alpha}\neq 0$. Note that if $d_{\alpha}\neq 0$ then $b_j\neq 0$.  Then $a_j={b_jc_{\alpha}\over d_{\alpha}}$. Taking into account condition (\ref{th1.1}), we have
\begin{equation*}
A_iB_j-B_iA_j = { A\over b_{r_i+1}} \cdot {b_jB} - {B\over a_{r_i+1}} \cdot {a_jA}= {AB ( a_{r_i+1}  b_{j} - b_{r_i+1} a_{j}) \over b_{r_i+1}  a_{r_i+1} } = {AB ( a_{r_i+1}  b_{j} - b_{r_i+1} {b_jc_{\alpha}\over d_{\alpha}}) \over b_{r_i+1}  a_{r_i+1} }= $$$$
= {AB b_j ( a_{r_i+1} d_{\alpha}   - b_{r_i+1} {c_{\alpha}}) \over d_{\alpha} b_{r_i+1}  a_{r_i+1} }
 \neq 0, \ i=\overline{0,k}, j=\overline{m+1,n}. 
\end{equation*}
If there exists a nonzero coefficient $A_{k+1}$ then condition (\ref{th1.10.2}) is fulfilled. Therefore $b_j\neq 0$ for $j=\overline{m+1,n}$.
\begin{equation*} 
	A_{k+1}B_j-B_{k+1}A_j = B_j=b_j B \neq 0, \quad j=\overline{m+1,n}.
\end{equation*}
If there exists a nonzero coefficient $B_{k+2}$ then condition (\ref{th1.10.1}) is fulfilled. Therefore $a_j\neq 0$ for $j=\overline{m+1,n}$.
\begin{equation*} 
A_{k+2}B_j-B_{k+2}A_j = - A_j=a_j A \neq 0, \quad j=\overline{m+1,n}.
\end{equation*}
Now we can solve equation (\ref{th1.17}).

Put $\rho_{i}=\mu_{i}*\overline{\mu}_{i}$, $\varrho_{j}=\mu_{\xi_j}*\overline{\mu}_{\xi_j}$. Then we have $\widehat\rho_{i}(y)=|\widehat\mu_{i}(y)|^2>0$, $\widehat\varrho_{j}(y)=|\widehat\mu_{\xi_j}(y)|^2>0$, $y\in Y$. The characteristic functions $\widehat\rho_{i}(y)$ and $\widehat\varrho_{j}(y)$ also satisfy equation (\ref{th1.17}), and all the factors in (\ref{th1.17}) are greater then zero.

Put $\varphi_i(y)= \log \widehat\rho_{i}(y)$, $\psi_j(y)= \log \widehat\varrho_{j}(y)$. It follows from (\ref{th1.17}) that the functions $\varphi_i(y)$ and $\psi_j(y)$ satisfy the following equation

\begin{equation}\label{th1.17.1}
\sum_{i=0}^{k+2}  \varphi_{i} \left(A_i u+B_i v\right) + \sum_{j=m+1}^{n}  \psi_{j} (A_j u+ B_j v) =\sum_{j=1}^{n}  \psi_{j} (C_j u+ D_j v),\quad u,v\in Y.
\end{equation}

It follows from Lemma \ref{FinDifEq} that for example the function $\varphi_0(y)$ satisfies the equation (\ref{fde2}) which takes the form

\begin{equation}\label{th1.17.2}
\Delta_{(A_0B_{1}-B_0A_{1})l_{1}} ... \Delta_{(A_0B_{k+2}-B_0A_{k+2})l_{m}} 
\Delta_{(A_0B_{m+1}-B_0A_{m+1})l_{1}} ... \Delta_{(A_0B_{n}-B_0A_{n})l_{m}}$$$$ \Delta_{(A_0D_{1}-B_0C_{1})k_{1}} ...  \Delta_{(A_0D_n-B_0C_n)k_n} \varphi_{0} (A_0 y)=0,\quad y\in Y,
\end{equation}
where $k_1, \dots, k_n$ and $l_{1}, \dots, l_{m}$ are arbitrary elements of $Y$. It follows from equation (\ref{th1.17.2}) and conditions (\ref{th1.10.-1}), (\ref{th1.13}), (\ref{th1.14}), (\ref{th1.18}), (\ref{th1.19}) that  

\begin{equation}\label{th1.17.3}
\Delta_{h}^{m+n+1} \varphi_{0} (y)=0,\quad y,h\in Y.
\end{equation}

Thus $\varphi_{0} (y)$ is a continuous polynomial on the group $Y$. Applying Lemma \ref{GrMarc}, we obtain that $\rho_0\in \Gamma(X)$. Hence by Lemma \ref{GrKramer} $\mu_0\in \Gamma(X)$ and respectively $\mu_{\xi_1}, \dots, \mu_{\xi_{r_1}} \in \Gamma(X)$. By analogy we obtain that all $\mu_{\xi_j}\in \Gamma(X)$ ($j=\overline{1,m}$).

(ii) Let $X=X_{(p)}$, where $p$ is a prime number ($p>2$). The proof of this case is the same as the proof of \textit{Case 1}(i). There is only one difference: if an integer $b_i
a_j - b_j a_i$ is not admissible for $X$ then $b_i a_j - b_j a_i=0\ (mod\ p)$.
At the end we note that since 
component of zero of the group $X$ is equal to zero, we have in this case
that $\Gamma(X)=D(X)$. 

\textit{Case 2.} Now we consider case 2 when all integers $b_i
a_j - b_j a_i$ for $i\neq j$ are not admissible for $X$.

(i) Let $X$ be a torsion-free group. In this case if all integers $b_i
a_j - b_j a_i$ are not admissible for $X$ then all $b_i a_j - b_j a_i=0$ for $i\neq j$.
Substitute $b_1y$ for $u$ and $-a_1y$ for $v$ in (\ref{th1.2}). Note that $a_jb_1 y - b_j a_1 y=0$ for all $j=\overline{1,n}$.  We get
\begin{equation}\label{th1.23}
1=\prod_{j=1}^{n}  \widehat\mu_{\xi_j} \left((b_1c_j-a_1d_j)y \right),\quad y\in Y.
\end{equation}
Taking into account condition (\ref{th1.1}),  we get from (\ref{th1.23}) that all $\widehat\mu_{\xi_j}(y)=1$ for all $y\in Y$, i.e. all $\mu_{\xi_j}\in D(X)$ ($j=\overline{1,n}$).

(ii) Let $X=X_{(p)}$, where $p$ is a prime number ($p>2$). In this case if all integers $b_i
a_j - b_j a_i$ are not admissible for $X$ then all $b_i a_j - b_j a_i=0\ (mod \ p)$ for $i\neq j$. Arguing as in \textit{Case 2} (i) we obtain that all $\mu_{\xi_j}\in D(X)$ ($j=\overline{1,n}$).
\hfill $\square$ 




\medskip

Theorem \ref{theorem1} is exact in the following sense. 

\begin{proposition}\label{theorem2}
	Let $X$ be a second countable locally
	compact Abelian group. If $X$ is not topologically isomorphic to any of the groups mentioned
	in Theorem $\ref{theorem1}$, then there exist independent random variables $\xi_j, j=\overline{1,n}, n \ge 2,$ with values in $X$ and distributions $\mu_{\xi_j}$
	with non-vanishing characteristic functions, 
	and integers
	$a_j, b_j, c_j, d_j $ such that the random vectors $(L_1,L_2)$ and $(L_3,L_4)$ are identically distributed, but for all those $i$, for which condition $(\ref{th1.1})$ is valid, $\mu_{\xi_i}\not\in\Gamma(X)$. The linear forms $L_1, L_2, L_3, L_4$ are defined by $(\ref{th1.1.1})-(\ref{th1.1.4})$.
\end{proposition}

\textit{\textbf{Proof.}} 
Let $x_0$ be an element of order $p$ in $X$. Let $K$ be a subgroup of $X$ generated by the element $x_0$. Then
$K\cong\mathbb{Z}(p)$.
Let $\xi_1$ and $\xi_2$ be independent
identically distributed random variables with values in $K$ and with
a distribution
\begin{equation}\label{th1.28}
\mu= m E_0+(1-m)E_{x_0},
\end{equation}
where $1/2 < m <1$. We consider the distribution $\mu$ as a distribution
on the group $X$. We have that
\begin{equation}\label{th1.29}
\widehat{\mu}(y)= m +(1-m) (x_0,y).
\end{equation}
Let $a_j=b_j=c_j=1$, $d_1=d_2=1-p$, i.e. $L_1=\xi_1+\xi_2$, $L_2=\xi_1+\xi_2$, $L_3=\xi_1+\xi_2$, $L_4=(1-p)\xi_1+(1-p)\xi_2$. It is easy to verify that condition (\ref{th1.1}) is valid for $i=1,2$. By Lemma \ref{Equation} the vectors $(L_1,L_2)$ and $(L_3,L_4)$ are identically distributed if
and only if the characteristic function $\widehat\mu(y)$ satisfy
equation (\ref{th1.2}) which takes the form
\begin{equation}\label{th1.28}
\widehat\mu(u+v)\widehat\mu(u+v)=\widehat\mu(u+(1-p)v)\widehat\mu(u+(1-p)v),
\quad u, v \in Y.
\end{equation}
Rewrite this equation in the following form
\begin{equation}\label{th1.30}
\widehat\mu(u+v)\widehat\mu(u+v)=\widehat\mu(u+v-pv)\widehat\mu(u+v-pv),
\quad u, v \in Y.
\end{equation}
Since $\sigma(\mu) \subset K \subset X_{(p)}$, we have $\widehat\mu(y+h)=\widehat\mu(y)$ for $y\in
Y^{(p)}$. Therefore equation (\ref{th1.30}) is an equality. Hence by Lemma \ref{Equation} the vectors $(L_1,L_2)$ and $(L_3,L_4)$ are identically distributed.
\hfill $\square$ 

\begin{remark}\label{cor1H}
	{\rm Let $L_1$, $L_2$ be as in Theorem \ref{theorem1}. Suppose that $L_3=L_1$ and $L_4=-L_2$. If the vectors $(L_1,L_2)$ and $(L_1,L_2)$ are identically distributed, then the conditional distribution of $L_2$ given $L_1$ is symmetric. Condition $(\ref{th1.1})$ takes the form
	\begin{equation}\label{cor1.1}
	a_i b_j + b_i a_j \text{ are admissible integers for } X \text{ for all } j=\overline{1,n}.
	\end{equation} 
	If we suppose that $a_j, b_j$ are admissible integers for $X$ then we obtain from  Theorem \ref{theorem1}
	the description of locally compact Abelian groups for which
	the group analogue of the Heyde theorem takes place (\cite[Theorem 3.1]{Myronyuk2020}).}
\end{remark}

\begin{remark}\label{cor2SD}
	{\rm We consider two sets of independent random variables $\xi_1, \dots, \xi_n$ and $\xi'_1, \dots, \xi'_n$, where $\xi_j$ and $\xi'_j$ are identically distributed.  Let $L_1$, $L_2$ be as in Theorem \ref{theorem1}. Suppose that $L_3=L_1$ and $L_4=L'_2=b_1\xi'_1+\cdots+b_n\xi'_n$. If the vectors $(L_1,L_2)$ and $(L_1,L'_2)$ are identically distributed, then $L_1$ and $L_2$ are independent. Condition $(\ref{th1.1})$ takes the form
		\begin{equation}\label{cor1.2}
		a_i b_i \text{ are admissible integers for all }i=\overline{1,n}.
		\end{equation} 
		We obtain from Theorem \ref{theorem1} 
		the description of locally compact Abelian groups for which
		the group analogue of the Darmois-Skitovich theorem takes place (\cite[Theorem 3]{Fe1990}).}
\end{remark}

\medskip

The following theorem describes locally compact Abelian groups for which the group
analogue of the Klebanov theorem takes place without assumption that
the characteristic functions of the considering distributions of random variables  do
not vanish.  Note that the obtained class of groups is changed
(compare with Theorem \ref{theorem1}).

\begin{theorem}\label{theorem3}
	Let $X=\mathbb{R}^n\times D$, where $n\geq 0$
	and $D$ is a discrete torsion-free group. Let $\xi_1, \dots, \xi_n$ be independent random variables with values in $X$ and distributions $\mu_{\xi_j}$. Consider the linear forms $(\ref{th1.1.1})$-$(\ref{th1.1.4})$,
	where the coefficients $a_j, b_j, c_j, d_j$ are integers. If the random vectors $(L_1,L_2)$ and $(L_3,L_4)$ are identically distributed, then for all those $i$ for which condition $(\ref{th1.1})$ is valid  $\mu_{\xi_i}\in\Gamma(X)$.	
\end{theorem}

To prove Theorem \ref{theorem3}, we need the following lemma.

\begin{lemma}\label{ConnectedCompact}
	Let $K$ be a connected compact Abelian
	group and $\mu_{j}$ be distributions on the group $K^*$. Assume that
	the characteristic functions $\widehat{\mu}_j(y)$  satisfy equation
	$(\ref{th1.2})$ on the group $K$  and $\widehat{\mu}_j(y)\geq 0$. Let
	$a_j, b_j, c_j, d_j$ be integers. For for all those $i$ for which $a_id_j-b_ic_j\neq 0$ for all $j=\overline{1,n}$, the characteristic function $\widehat\mu_i(y)=1$, $y\in K$.
\end{lemma}

The proof of Lemma \ref{ConnectedCompact}  is carried out in the same way as the proof of Lemma 3.1 of the paper \cite{Fe1990} or the proof of Lemma 3.8 of the paper \cite{Myronyuk2020}, but in our case it is based on the Klebanov
theorem on the real line. We omitted the prove of Lemma \ref{ConnectedCompact} in our paper.

\medskip

\textit{\textbf{Proof of Theorem \ref{theorem3}.}} Lemma \ref{Equation} implies that the characteristic functions of distributions
$\mu_{\xi_j}$ satisfy equation (\ref{th1.2}). Put $\nu_{j}=\mu_{\xi_j}*\overline{\mu}_{\xi_j}$. Then we have $\widehat\nu_{j}(y)=|\widehat\mu_{\xi_j}(y)|^2\geq 0 $, $y\in Y$. The characteristic functions $\widehat\nu_{j}(y)$ also satisfy equation (\ref{th1.2}). If we prove that $\nu_j\in \Gamma(X)$, then it follows from Lemma \ref{GrKramer}, that $\mu_{\xi_j}\in \Gamma(X)$. Therefore we can assume from the beginning that all factors $\widehat\mu_{\xi_j}(y)\geq 0$  in (\ref{th1.2}).

Renumbering functions in (\ref{th1.2}), we can assume that condition (\ref{th1.1}) is fulfilled for all $i=\overline{1,m}$ and it is not fulfilled for all $i=\overline{m+1,n}$, $m\leq n$.

Since $X= \mathbb{R}^n\times D$, we have $Y\approx \mathbb{R}^n\times K$, where $K$ is a connected
compact group.  We can assume without loss of generality that
$Y=\mathbb{R}^n\times K$.

Consider the restriction of equation (\ref{th1.2}) to $K$. It follows
from Lemma \ref{ConnectedCompact} that $\widehat\mu_{\xi_i}(y)=1$ on $K$ for $i=\overline{1,m}$. Then $\widehat\mu_{\xi_i}(y+h)=\widehat\mu_{\xi_i}(y)$ for all $y\in Y$, $h\in K$. We consider the
restriction of equation (\ref{th1.2}) to each one-dimensional subspace in $\mathbb{R}^n$
and obtain from the Klebanov theorem that all restriction of the
characteristic functions $\widehat\mu_{\xi_i}(y)$ are the characteristic
functions of Gaussian distributions for $i=\overline{1,m}$. Therefore $\widehat\mu_{\xi_i}(y)$ are the characteristic
functions of Gaussian distributions for $y\in \mathbb{R}^n$ and $i=\overline{1,m}$. Hence
$\mu_{\xi_i}\in\Gamma(X)$ for $i=\overline{1,m}$. \hfill $\square$ 

\medskip

Theorem \ref{theorem3} is exact in the following sense. 

\begin{proposition}\label{theorem4}
	Let $X$ be a second countable locally
	compact Abelian group. If $X$ is not topologically isomorphic to the group mentioned
	in Theorem $\ref{theorem3}$, then there exist independent random variables $\xi_j, j=\overline{1,n}, n \ge 2,$ with values in $X$ and distributions $\mu_{\xi_j}$, and integers
	$a_j, b_j, c_j, d_j $ such that the random vectors $(L_1,L_2)$ and $(L_3,L_4)$ are identically distributed, but for all those $i$, for which condition $(\ref{th1.1})$ is valid, $\mu_{\xi_i}\not\in\Gamma(X)*I(X)$. The linear forms $L_1, L_2, L_3, L_4$ are defined by $(\ref{th1.1.1})-(\ref{th1.1.4})$.
\end{proposition}

To prove proposition \ref{theorem4} we need the following lemma and remark.

\begin{lemma}\label{HeydeGr} $(\cite{Myronyuk2020})$
	Let $X$ be a second countable locally
	compact Abelian group. If either $X\not\approx\mathbb{R}^n\times D$, where $n\geq 0$
	and $D$ is a discrete torsion-free group, or  $X\neq X_{(2)}$, or $X\neq X_{(3)}$, then there exist independent random variables $\xi_j, j=\overline{1,n}, n \ge 2,$ with values in $X$ and distributions $\mu_j$,
	and admissible integers $a_j, \ b_j $ such that $b_i a_j + b_j a_i$
	are admissible integers for $X$ for all $i, j$, such that the
	conditional distribution of $L_2$ given
	$L_1$ is symmetric, but all
	$\mu_j\not\in\Gamma(X)*I(X)$. The linear forms $L_1, L_2$ are defined by $(\ref{th1.1.1})-(\ref{th1.1.2})$.
\end{lemma}

\begin{lemma}\label{KlebanovZ2} 
	Let $X=\mathbb{Z}(2)$ or $X=\mathbb{Z}(3)$. Then there exist independent random variables $\xi_j, j=\overline{1,n}, n \ge 2,$ with values in $X$ and distributions $\mu_{\xi_j}$, and integers
	$a_j, b_j, c_j, d_j $ such that the random vectors $(L_1,L_2)$ and $(L_3,L_4)$ are identically distributed, but for all those $i$, for which condition $(\ref{th1.1})$ is valid, $\mu_i\not\in I(X)$. 
\end{lemma}

\textit{\textbf{Proof.}} Let $L_1=\xi_1+...+\xi_{n-2}+\xi_{n-1}$, $L_2=\xi_1+...+\xi_{n-2}+\xi_{n}$, $L_3=\xi_1+...+\xi_{n-2}+\xi_{n-1}$, $L_4=\xi_{n}$. Condition $(\ref{th1.1})$ is valid for all $i=\overline{1,n-2}$. Indeed, we have in this case $a_id_j-b_ic_j=d_j-c_j=\pm 1$. Note that condition $(\ref{th1.1})$ is not valid for $i=n-1$ and $i=n$. Suppose that $\mu_{\xi_1}, \dots, \mu_{\xi_{n-2}}$ are arbitrary distributions. Put $\mu_{\xi_{n-1}}=\mu_{\xi_{n}}=m_X$  and verify that the random vectors $(L_1,L_2)$ and $(L_3,L_4)$ are identically distributed. By Lemma 1 it
suffices to show that the characteristic functions of distributions $\mu_{\xi_{j}}$ satisfy equation
(\ref{th1.2}) which takes the form
\begin{equation}\label{lemz2.1}
	\prod_{j=1}^{n-2}  \widehat\mu_{\xi_j} (u+ v) \widehat{m}_X(u)\widehat{m}_X(v)=\prod_{j=1}^{n-2}  \widehat\mu_{\xi_j} (u) \widehat{m}_X(u)\widehat{m}_X(v),\quad u,v\in Y.
\end{equation}
Obviously, equation (\ref{lemz2.1}) is satisfied if $u=v=0$. Suppose that either $u\neq 0$ or $v\neq 0$. Then $\widehat{m}_X(u)\widehat{m}_X(v)=0$, i.e. the
left-hand side and the right-hand side of equation (\ref{lemz2.1}) is equal to zero. It is clear that we can choose distributions $\mu_{\xi_1}, \dots, \mu_{\xi_{n-2}}$ in such a way that $\mu_{\xi_j}\notin I(X)$.
\hfill $\square$ 

\medskip

\textbf{\textit{Proof of Proposion \ref{theorem4}.} } Suppose that either $X\not\approx\mathbb{R}^n\times D$, where $n\geq 0$
and $D$ is a discrete torsion-free group, or  $X\neq X_{(2)}$, or $X\neq X_{(3)}$. Let $L_1$, $L_2$ be as in (\ref{th1.1.1}) and (\ref{th1.1.2}). We set $L_3=L_1$, $L_4=-L_2$. 
Then condition $(\ref{th1.1})$ takes the form (\ref{cor1.1}).
The random vectors $(L_1, L_2)$ and $(L_1, - L_2)$ are identically distributed if and only if the conditional distribution of $L_2$ given $L_1$ is symmetric. Now the statement of Proposition \ref{theorem4} follows from Lemma \ref{HeydeGr}.

Suppose now that $X= X_{(2)}$ or $X= X_{(3)}$. Then the statement of Proposition \ref{theorem4} follows from Lemma \ref{KlebanovZ2}.
\hfill $\square$ 

\begin{remark}\label{cor3H}
	{\rm As in Remark \ref{cor1H}  we obtain from  Theorem \ref{theorem3} that 
		the group analogue of the Heyde theorem takes place for groups $X=\mathbb{R}^n\times D$, where $n\geq 0$
		and $D$ is a discrete torsion-free group (compare with \cite[Theorem 3.6]{Myronyuk2020}). But it was proved in \cite{Myronyuk2020} that the group analogue of the Heyde theorem takes place for groups $X=X_{(3)}$. }
		
	{\rm In general the group analogue of the Klebanov theorem does not take place for such groups. Let us suppose that $X=X_{(3)}$ and in addition all $a_j, b_j, c_j, d_j$ are admissible integers for $X$. Then we can assume without loss of generality
		that $a_j, b_j, c_j, d_j$ are equal to $\pm 1$. Suppose also that condition $(\ref{th1.1})$ is fulfilled for $i=\overline{1,m}$, where $m\leq n$.  It follows from condition $(\ref{th1.1})$ that if $a_i=-b_i$ ($1\leq i \leq m$) then $c_i= d_i$ for all $i=\overline{1,n}$ or if $a_i=b_i$ then $c_i= -d_i$ for all $i=\overline{1,n}$. For definiteness, suppose that $a_i=-b_i \text{ for all } i=\overline{1,m}$. Then it follows from $(\ref{th1.1})$ that $c_i= d_i \text{ for all } i=\overline{1,n}$. Thus equation (\ref{th1.2}) has the form
		\begin{equation}\label{th3.1.cor}
		\prod_{j=1}^{m}  \widehat\mu_{\xi_j} (a_j (u-v))\prod_{j=m+1}^{n}  \widehat\mu_{\xi_j} (a_j (u+v))=\prod_{j=1}^{n}  \widehat\mu_{\xi_j} (c_j (u+v)),\quad u,v\in Y.
		\end{equation}
		Putting $v=-u$ in (\ref{th3.1.cor}), we get
		\begin{equation}\label{th3.2.cor}
		\prod_{j=1}^{m}  \widehat\mu_{\xi_j} (2a_j u)=1,\quad u\in Y.
		\end{equation}
		Since $f_{2a_j}\in Aut(Y)$ we get from (\ref{th3.2.cor}) that all $\widehat\mu_{\xi_j} (y)=1$ for $y\in Y$. Therefore all $\mu_{\xi_j}\in D(X)$. Thus, in this special case the group analogue of the Klebanov theorem implies  the group analogue of the Heyde theorem for $X=X_{(3)}$.}
\end{remark}

\begin{remark}\label{cor4SD}
	{\rm As in Remark \ref{cor2SD}  we obtain from  Theorem \ref{theorem3} that 
		the group analogue of the Darmois-Skitovich theorem takes place for groups $X=\mathbb{R}^n\times D$, where $n\geq 0$
		and $D$ is a discrete torsion-free group (compare with \cite[Theorem 1]{Fe1990}). But it was proved in \cite{Fe1990} that the group analogue of the Darmois-Skitovich theorem takes place for groups $X=X_{(2)}$. }
	
	{\rm In general the group analogue of the Klebanov theorem does not take place for such groups. Let us suppose that $X=X_{(2)}$ in Theorem \ref{theorem3}. In addition we suppose that $n=2k$ and $a_j, b_j, c_j$ are admissible for only $j=\overline{1,k}$ and $d_j$ are admissible for only $j=\overline{k+1,n}$. Then we can assume without loss of generality
		that $a_j=b_j=c_j=1$ for $j=\overline{1,k}$, $a_j=b_j=c_j=0$ for $j=\overline{k+1,n}$, and $d_j=1$ for $j=\overline{k+1,n}$, $d_j=0$ for $j=\overline{1,k}$.
		Thus equation (\ref{th1.2}) has the form
		\begin{equation}\label{th3.1.cor4}
		\prod_{j=1}^{k}  \widehat\mu_{\xi_j} (u+v)=\prod_{j=1}^{k}  \widehat\mu_{\xi_j} (u) \prod_{j=k+1}^{n}  \widehat\mu_{\xi_j} (v),\quad u,v\in Y.
		\end{equation}
		Putting $v=u$ in (\ref{th3.1.cor4}), we get
		\begin{equation}\label{th3.2.cor4}
		\prod_{j=1}^{n}  \widehat\mu_{\xi_j} (u)=1,\quad u\in Y.
		\end{equation}
		We get from (\ref{th3.2.cor4}) that all $\widehat\mu_{\xi_j} (y)=1$ for $y\in Y$. Therefore all $\mu_{\xi_j}\in D(X)$. Thus, in this special case the group analogue of the Klebanov theorem implies  the group analogue of the Darmois-Skitovich theorem for $X=X_{(2)}$.}
\end{remark}

\begin{remark}
	{\rm Theorems \ref{theorem1} and \ref{theorem3} fail if we omit condition $(\ref{th1.1})$. Indeed, let $\xi_1$ and $\xi_2$ be independent identically distributed random variables with a distribution $\mu$. Put $L_1=L_3=\xi_1+\xi_2$, $L_2=-L_4=\xi_1-\xi_2$. In this case condition $(\ref{th1.1})$ is not fulfilled.  If the random vectors $(L_1, L_2)$ and $(L_1, - L_2)$ are identically distributed then Lemma \ref{Equation} implies that the characteristic function $\widehat{\mu}(y)$ satisfies equation (\ref{th1.2}) which takes the form}
	\begin{equation}\label{rem5.1}
	\widehat\mu(u+v)\widehat\mu(u-v)=\widehat\mu(u-v)\widehat\mu(u+v),
	\quad u, v \in Y.
	\end{equation}
	{\rm Equation (\ref{rem5.1}) is the identity for any function $\widehat\mu(y)$. Therefore we can suppose that $\mu\not\in\Gamma(X)*I(X)$.	}
\end{remark}

\section{Case of Q-independent random variables}

In the article  \cite{KS} A.M.Kagan and G.J.Sz\'ekely introduced a
notion of $Q$-independence of random variables which generalizes the
notion of independence of random variables. Then in \cite{Fe2017}
G.M. Feldman in a natural way introduced a notion of
$Q$-independence of random variables taking values in a locally
compact Abelian group. These studies
were continued in \cite{My2019} and in \cite{Myronyuk2020}.

Let $\xi_1, \dots, \xi_n$ be random variables with values in the
group $X$. The random variables $\xi_1, \dots, \xi_n$
are $Q$-independent if the characteristic function of the vector
$(\xi_1, \dots, \xi_n)$ can be represented in the form
\begin{equation}\label{s2.1}
\widehat\mu_{(\xi_1, \dots, \xi_n)}(y_1, \dots, y_n)=\prod_{j=1}^n\widehat\mu_{\xi_j}(y_j)\exp\{q(y_1,
\dots, y_n)\}, \quad y_j\in Y,
\end{equation}
where $q(y_1, \dots, y_n)$ is a continuous polynomial on the group
$Y^n$ such that $q(0, \dots, 0)=0$.

We formulate an analogue of Lemma \ref{Equation} for $Q$-independent random variables. 

\begin{lemma}\label{EquationQ}
	Let $X$ be a second countable locally
	compact Abelian group, $Y=X^*$. Let $\xi_1, \dots, \xi_n$ be $Q$-independent random variables with values in $X$ and distributions $\mu_{\xi_j}$. Consider the linear forms $(\ref{th1.1.1})-(\ref{th1.1.4})$
	where the coefficients $a_j, b_j, c_j, d_j$ are integers. The random vectors $(L_1,L_2)$ and $(L_3,L_4)$ are identically distributed if and only if the characteristic functions $\widehat{\mu}_{\xi_j}(y)$ satisfy the equation
	\begin{equation}\label{th1.2Q}
	\prod_{j=1}^{n}  \widehat\mu_{\xi_j} (a_j u+ b_j v)=\prod_{j=1}^{n}  \widehat\mu_{\xi_j} (c_j u+ d_j v) \exp\{q(u,v)\},\quad u,v\in Y,
	\end{equation}
	where $q(u,v)$ is a continuous polynomial on the group
	$Y^2$ such that $q(0, 0)=0$.
\end{lemma}

\textit{\textbf{Proof.}} Using the $Q$-independence of $\xi_1, \dots, \xi_n$, in the same way as in Lemma \ref{Equation}, we obtain from the definition that
\begin{eqnarray*}
	\widehat\mu_{(L_1,L_2)}(u,v)=\mathbf{E}\left[ (\xi_1,a_1 u +b_1 v)\cdots(\xi_n, a_n u+b_n v) \right]= \\ = \prod_{j=1}^{n}  \widehat\mu_{\xi_j} (a_j u+ b_j v) \exp\{q_1(a_1 u +b_1 v,\dots, a_n u +b_n v)\},\quad u,v\in Y.
\end{eqnarray*}
Analogously we obtain that
\begin{equation*}
\widehat\mu_{(L_3,L_4)}(u,v)=\prod_{j=1}^{n}  \widehat\mu_{\xi_j} (c_j u+ d_j v) \exp\{q_2(c_1 u+ d_1 v,\dots, c_n u+ d_n v)\},\quad u,v\in Y.
\end{equation*} 
Since the vectors $(L_1,L_2)$ and $(L_3,L_4)$ are identically distributed, we have $\widehat\mu_{(L_1,L_2)}(u,v)=\widehat\mu_{(L_3,L_4)}(u,v)$. Put $q(u,v)=q_2(c_1 u+ d_1 v,\dots, c_n u+ d_n v)-q_1(a_1 u +b_1 v,\dots, a_n u +b_n v)$. Thus, equation
(\ref{th1.2Q}) is valid.
\hfill $\square$

\medskip

Theorem \ref{theorem1} holds true if we consider $Q$-independent random variables instead of independent random variables. 

\begin{theorem}\label{theorem1Q}
	
	Let $X$ be a second countable locally
	compact Abelian group. Let $\xi_1, \dots, \xi_n$ be $Q$-independent random variables with values in $X$ and distributions $\mu_{\xi_j}$ with non-vanishing characteristic functions. Consider the linear forms $(\ref{th1.1.1})-(\ref{th1.1.4})$, 
	where the coefficients $a_j, b_j, c_j, d_j$ are 
	integers. If the random vectors $(L_1,L_2)$ and $(L_3,L_4)$ are identically distributed, then for all those $i$, for which condition $(\ref{th1.1})$ is valid,
	the following
	statements hold:
	
	\qquad $(i)$ If $X$ is a torsion-free group then 
	$\mu_{\xi_i}\in\Gamma(X)$;
	
	\qquad $(ii)$ If $X=X_{(p)}$ then $\mu_{\xi_i}\in D(X)$.
\end{theorem}

\textit{\textbf{Proof.}} The proof of (i) of Theorem \ref{theorem1Q} is almost the same as the proof of Theorem \ref{theorem1}. Instead of equation (\ref{th1.2}) we have to consider equation (\ref{th1.2Q}). Instead of equation (\ref{th1.17.1}) we get the following equation

\begin{equation}\label{th1.17.1Q}
\sum_{i=0}^{m}  \varphi_{i} \left(A_i u+B_i v\right) =\sum_{j=1}^{n}  \psi_{j} (C_j u+ D_j v)+q(u,v),\quad u,v\in Y.
\end{equation}
Then we apply Lemma \ref{FinDifEq} and complete the proof of (i) of Theorem \ref{theorem1Q} in the same way as Theorem \ref{theorem1}.

The proof of (ii) of Theorem \ref{theorem1Q} follows from Theorem \ref{theorem1} and the following fact. It is well known that any continuous polynomial is equal to the constant
on compact elements (see e.g. \cite[\S5.7]{FeBook2}). If the group
$Y$ consists only of compact elements then the connected component
of zero of the group $ X $ is equal to zero (see e.g.
\cite[\S24.17]{HeRo1}). Thus, independence and Q-independence of
random variables are equivalent on groups $X$ whose connected
component of zero is equal to zero. Hence statement of (ii) of Theorem \ref{theorem1Q} follows from Theorem \ref{theorem1}. 
\hfill $\square$

\medskip

Theorem \ref{theorem3} also holds true if we consider $Q$-independent random variables instead of independent random variables. 

\begin{theorem}\label{theorem3Q}
	Let $X=\mathbb{R}^n\times D$, where $n\geq 0$
	and $D$ is a discrete torsion-free group. Let $\xi_1, \dots, \xi_n$ be $Q$-independent random variables with values in $X$ and distributions $\mu_{\xi_j}$. Consider the linear forms $(\ref{th1.1.1})$-$(\ref{th1.1.4})$,
	where the coefficients $a_j, b_j, c_j, d_j$ are integers. If the random vectors $(L_1,L_2)$ and $(L_3,L_4)$ are identically distributed, then for all those $i$ for which condition $(\ref{th1.1})$ is valid  $\mu_{\xi_i}\in\Gamma(X)$.	
\end{theorem}

\begin{remark}
	{\rm Since independent random variables are also $Q$-independent, Propositions \ref{theorem2} and \ref{theorem4} prove the exactness of Theorems \ref{theorem1Q} and \ref{theorem3Q}.}
\end{remark}

\section*{Acknowledgements}

The author would like to thank the Volkswagen Foundation (VolkswagenStiftung), the Bielefeld University and Prof. Dr. Friedrich G\"{o}tze for the support and warm reception.

\section*{Funding}

This research was supported by VolkswagenStiftung - Az. 9C108.

\end{document}